\documentclass[a4paper,11pt,twoside]{amsart}
\usepackage[latin1]{inputenc}
\usepackage[T1]{fontenc}
\usepackage{amsmath}
\usepackage{amsthm}
\usepackage{amssymb}
\usepackage[all]{xy}

\newtheorem{thm}{Theorem}[section]

\newtheorem{prop}[thm]{Proposition}

\newtheorem{lem}[thm]{Lemma}
\newtheorem{cor}[thm]{Corollary}

\numberwithin{equation}{section}

\theoremstyle{definition}
\newtheorem{definition}[thm]{Definition}
\newtheorem{remark}[thm]{Remark}

\newcommand{\qqed}{\hspace*{\fill}$\Box$}
\newcommand{\supp}{\operatorname{Supp}}

\newcommand{\Db}{{\rm D}^{\rm b}}

\newcommand{\TAut}{{\rm Aut}^{\rm t}}
\newcommand{\Aut}{{\rm Aut}}

\newcommand{\Gl}{{\rm Gl}}

\newcommand{\NS}{{\rm NS}}

\newcommand{\coh}{{\rm Coh}}

\newcommand{\End}{{\rm End}}
\newcommand{\Hom}{{\rm Hom}}

\newcommand{\ext}{{\rm ext}}

\newcommand{\Stab}{{\rm Stab}}

\newcommand{\id}{{\rm id}}

\newcommand{\Ext}{{\rm Ext}}

\newcommand{\cal}{\mathcal}
\newcommand{\ka}{{\cal A}}

\newcommand{\kd}{{\cal D}}

\newcommand{\kh}{{\cal H}}

\newcommand{\ko}{{\cal O}}
\newcommand{\kp}{{\cal P}}

\newcommand{\kt}{{\cal T}}

\newcommand{\ks}{{\cal S}}

\newcommand{\ZZ}{\mathbb{Z}}
\newcommand{\QQ}{\mathbb{Q}}
\newcommand{\RR}{\mathbb{R}}
\newcommand{\CC}{\mathbb{C}}

\newcommand{\PP}{\mathbb{P}}

\renewcommand{\to}{\xymatrix@1@=15pt{\ar[r]&}}
\renewcommand{\rightarrow}{\xymatrix@1@=15pt{\ar[r]&}}
\renewcommand{\mapsto}{\xymatrix@1@=15pt{\ar@{|->}[r]&}}
\renewcommand{\twoheadrightarrow}{\xymatrix@1@=15pt{\ar@{->>}[r]&}}
\renewcommand{\hookrightarrow}{\xymatrix@1@=15pt{\ar@{^(->}[r]&}}
\newcommand{\congpf}{\xymatrix@1@=15pt{\ar[r]^-\sim&}}
\renewcommand{\cong}{\simeq}

\begin{document}
\title{Stability conditions via spherical objects}
\author{D. Huybrechts}
\address{Mathematisches Institut,
Universit{\"a}t Bonn, Endenicher Alle 60, 53115 Bonn, Germany}
\email{huybrech@math.uni-bonn.de} 
\thanks{This work was supported by the SFB/TR 45 `Periods,
Moduli Spaces and Arithmetic of Algebraic Varieties' of the DFG
(German Research Foundation). The hospitality and financial support
of the Mathematical Institute Oxford is gratefully acknowledged.
}
\maketitle

\begin{abstract}
An object in the bounded derived category $\Db(X)$ of coherent sheaves on a complex projective K3 surface $X$ is spherical if it is rigid and simple. Although spherical objects form only a discrete set in the moduli stack
of complexes, they determine much of the structure of $X$ and $\Db(X)$. Here we show
that a stability condition  on $\Db(X)$ is determined by the stability of spherical objects. 
\end{abstract}

\medskip

Consider the bounded derived category $\Db(X)$ of the abelian category $\coh(X)$ of coherent sheaves on a complex projective K3 surface $X$.
An object $A$ in $\Db(X)$, i.e.\ a bounded complex of coherent sheaves on $X$, is called \emph{spherical} if $$\Ext^*(A,A)\cong H^*(S^2,\CC).$$

Spherical objects play a central role in the theory of K3 surfaces. In the classical theory they occur as
$(-2)$-curves, i.e.\ smooth rational curves $\PP^1\cong C\subset X$, and indeed
the structure sheaf $\ko_C$ of a $(-2)$-curve
$C$ is  a spherical object in $\Db(X)$. Other examples of spherical objects in $\Db(X)$ are provided by
 line bundles $L$ or, more generally, rigid stable vector bundles.


More recently, spherical objects $A$ and their associated spherical twists $T_A$ have
been used to give a conjectural description of
the group of all exact linear autoequivalences $\Aut(\Db(X))$ (see \cite{BK3}). In this language, the reflection $s_C$ classically associated
to a $(-2)$-curve $C$ and used to generate the Weyl group $W_X$ of a K3 surface, can be reinterpreted as the cohomological
action of the spherical twist $T_{\ko_C(-1)}$. Spherical objects also seem to play
a role in the study of Chow groups and  the arithmetic of K3 surfaces (see \cite{HCH,HuyBB}).

Clearly,  point sheaves $k(x)$, which are semirigid  objects in $\Db(X)$, determine the geometry of $X$ completely. But it seems that the much smaller discrete
set $\ks\subset\Db(X)$ of all spherical objects carries essentially the same information.

\medskip

The purpose of the present article is to stress this point  further by showing that a stability condition
on the derived category $\Db(X)$ is determined by the phases of only
spherical objects. 

Stability conditions as introduced by Bridgeland in \cite{BSt} have been studied intensively
for  K3 categories. In \cite{BK3} Bridgeland studies a distinguished connected
component $\Sigma$ of the space $\Stab(X):=\Stab(\Db(X))$ (as before, $X$ a projective K3 surface) of all stability conditions and 
conjectures that $\Sigma$ is simply-connected and preserved by the group of exact autoequivalences
$\Aut(\Db(X))$. For generic non-projective K3 surfaces $\Stab(X)$ was completely described
in \cite{HMS}. 

The case of local K3 surfaces is more accessible. More precisely,
for the minimal resolution $\pi:X\to\CC^2/G$ of a Kleinian singularity one can consider K3
categories $\kd\subset\hat\kd\subset\Db(X)$ of complexes supported on the exceptional
divisor (resp.\ with vanishing $R\pi_*$). The spaces $\Stab(\kd)$ and $\Stab(\hat\kd)$ 
are studied in detail in \cite{Ishi,Th} for  $A_n$-singularities and in \cite{Brav,BKl} in general.
Roughly, the analogue of Bridgeland's conjecture, originally  formulated for
projective K3 surfaces, are known to hold in the local situation.


There are at least two reasons why stability conditions in the local situation 
can be understood almost completely while the case of projective K3 surfaces still eludes
us. Firstly, the group of autoequivalences of $\Db(X)$ for a projective K3 surface $X$ is
much more complex than `just' a braid group. Indeed, $\Db(X)$ can host many different $A_n$-configurations of spherical objects at a time, which might be interlinked in a complicated manner.
Secondly, in the local case the categories under consideration are generated by spherical objects
and, in particular, their Grothendieck groups are of finite rank. A priori, the structure of $\Db(X)$
for a projective K3 surface seems more complicated due to the many
objects not generated by spherical objects.

The goal of this paper is to show that also for a projective K3 surface $X$ the space of stability 
conditions on $\Db(X)$ can be studied purely in terms of a configuration of spherical objects,
in other words in terms of a category that is spanned by spherical objects. In some sense this
is meant to bridge the gap  between the existing work in the local and in the global
setting, but whether it can be useful to prove Bridgeland's conjecture remains to be seen.

\medskip

The main result of the paper (see Theorem \ref{thm:main}) is concerned with two stability conditions  $\sigma=(\kp,Z)$
and $\sigma'=(\kp',Z')$ in the distinguished connected component
$\Sigma$ of the space of all stability conditions $\Stab(X)$.

\begin{thm}\label{thmMAIN}
Assume $Z=Z'$. Then $$\sigma=\sigma'$$
if and only if for all spherical objects $A$: 
$$A \text{ is }\sigma\text{-stable of
phase }\varphi\text{ if and only if }A\text{ is }\sigma'\text{-stable of phase }\varphi.$$ 
\end{thm}

The result can be reformulated in terms of a new metric on $\Stab(X)$, only taking spherical objects into account, which  
by the theorem turns out to be equivalent to the one defined by Bridgeland in \cite{BSt} (see Corollary \ref{cor:spmetr}).

This point of view is the motivation for the following construction. Consider the full triangulated
subcategory $\ks^*\subset\Db(X)$ generated by $\ks$. Note that in generating $\ks^*$ we do not allow taking
direct summands (see \cite{HuyBB} for details). Then $\ks^*$ is dense in $\Db(X)$ and its Grothendieck
group $K(\ks^*)\subset K(X)$ equals $N(X)=\ZZ\oplus \NS(X)\oplus \ZZ$ (under the additional
but presumably superfluous assumption $\rho(X)\geq2$).

  The triangulated category $\ks^*$ does not carry a bounded
t-structure and, therefore, no stability condition. But considering a weaker notion of stability conditions one
can introduce $\Stab(\ks)=\Stab(\ks^*)$ which as in \cite{BSt} is endowed with a natural (generalized) metric
(see Section \ref{sect:bfks}). The restriction
of a stability condition on $\Db(X)$ to $\ks^*\subset \Db(X)$ is then well-defined, i.e.\ there exists
a continuous map $\Stab(X)\to\Stab(\ks)$. 

As a consequence of Theorem \ref{thmMAIN} one obtains the following   (see Corollary \ref{cor:stabsph})
\begin{cor}
On the distinguished component the restriction yields an embedding
$$\Sigma\,\,\hookrightarrow \Stab(\ks)$$
which identifies the natural metric on $\Stab(\ks)$ with the spherical metric
$d_\ks$ on $\Sigma$.
\end{cor}

Note that for $\ks^*$ there is no difference between the Grothendieck group of $\ks^*$ and its quotient
by numerical equivalence $\sim_\ks$. Thus $K(\ks^*)/_{\sim_\ks}\cong K(X)/_\sim=N(X)$ and, therefore, maximal components
of $\Stab(\ks)$ and $\Stab(X)$ are modeled locally over the same linear space.

\medskip

Here is an outine of the paper.
Section \ref{sect:General} contains the basic definitions and results on stability conditions on $\Db(X)$
and explains some useful techniques  (Lemma \ref{lem:usefulsheaves}, \ref{lem:Lemma2})
to study the heart of a standard stability condition. In Section \ref{sect:Spherical} we recall that stable factors
of a spherical objects are again spherical and study spherical objects in the heart of a standard
stability condition. Section \ref{sect:Main} contains the proof of the main theorem. It is first proved for
the case that one of the stability conditions is standard. The generic case can be reduced to this
by applying autoequivalences, but the case of stability conditions in the boundary of the set of standard stability conditions is more complicated. The result can be rephrased in terms of the spherical metric,
which  is explained in Section \ref{sect:Sphmetr}. The last part could be read together with \cite{HuyBB} which discusses the category $\ks^*$ from a different angle and in more detail. The appendix collects a few observations
on the groups $\Aut(\Db(X))$ and $\Aut(X)$.


\section{General remarks on stability conditions}\label{sect:General}

\subsection{}\label{sect:basics} Recall that a stability condition $\sigma=(\kp,Z)$ on a triangulated category $\kd$ 
as introduced by Bridgeland in \cite{BSt} consists of a slicing $\kp$ and a stability function $Z$. 

The \emph{slicing} $\kp$ of $\sigma$ is given by  full abelian subcategories $\kp(\phi)\subset\kd$, $\phi\in\RR$.
The slicing has  two properties:   $$\Hom(\kp(\phi_1),\kp(\phi_2))=0\text{~ for~} \phi_1>\phi_2\text{~and~} \kp(\phi)[1]=\kp(\phi+1).$$ 
Objects in $\kp(\phi)$ are the \emph{semistable objects of phase $\phi$}. 
Let $\kp(\phi)^s\subset\kp(\phi)$ denote the subcategory of all \emph{stable} objects $E\in\kp(\phi)$, i.e.\ objects
$E\in\kp(\phi)$ not containing any proper non-trivial subobject in $\kp(\phi)$. 

The \emph{stability function}  is a linear function $Z:K(\kd)\to\CC$ such that $Z(E)\in\exp(i\pi\phi)\RR_{>0}$ for
all $0\ne E\in\kp(\phi)$.

\medskip

We shall only consider locally finite numerical stability
conditions. The latter means that the stability function
factors as  $Z:K(\kd)\to\Lambda\to\CC$ with $\Lambda:=K(\kd)/_\sim$ the numerical
Grothendieck group of $\kd$. 

The finiteness of the stability condition is a technical assumption that in \cite{BSt} enters  the discussion of the topology on the
space of stability conditions. Here,  finiteness will explicitly only be used to ensure
that the abelian  categories $\kp(\phi)$ are of finite length, i.e.\ any semistable
object has a finite filtrations with stable quotients.

In particular, for any stability condition $\sigma=(\kp,Z)$ and any object $0\ne E\in\kd$ there exists a \emph{$\sigma$-stable decomposition}, i.e.\  a diagram of exact triangles 
\begin{equation}\label{eqn:stabledomp}
\xymatrix{&F_1\ar@{=}[dl]\ar[rr]&&F_2\ar[dl]\ar[r]&\ldots\ar[r]& F_{m-1}\ar[rr]&&F_m\ar@{=}[r]\ar[dl]&E\\
A_1&&A_2\ar[ul]^{[1]}&&&&A_m\ar[ul]^{[1]}&&&}\end{equation}
with $A_i\in\kp(\phi_i)$ and such that $\phi_1\geq\ldots\geq\phi_m$. The minimal and maximal phases of $E$ are defined
as $\phi^-(E):=\phi_m$ resp.\ $\phi^+(E):=\phi_1$;  
they are uniquely determined. The $A_i$ are called the $\sigma$-stable factors of $E$ and they are unique up to permutation among those of the same phase. Note that as in the classical case the two morphisms
$A_1\to E$ and $E\to A_m$ are always non-trivial.

Requiring strict inequalities leads to the Harder--Narasimhan (or $\sigma$-semistable) decomposition by $\sigma$-semistable factors. This decomposition is unique. 

\subsection{}\label{subsect:stabMet} Let $\Stab(\kd)$ be the space of all (locally finite and numerical)
stability conditions on $\kd$. In \cite{BSt} Bridgeland uses a generalized metric to define a topology on $\Stab(\kd)$.
The distance between two slicings $\kp$ and $\kp'$ is measured by
$$f(\kp,\kp'):=\sup_{0\ne E\in\kd}\{|\phi^+(E)-{\phi'}^+(E)|,|\phi^-(E)-{\phi'}^-(E)|\}$$
where $\phi^\pm$ and ${\phi'}^\pm$ are the minimal (resp.\ maximal) phases with respect to $\kp$ resp.\ $\kp'$.

The generalized metric $d(\sigma,\sigma')$ between two stability conditions $\sigma=(\kp,Z),\sigma'=(\kp',Z')\in \Stab(X)$
combines $f(\kp,\kp')$ with a distance function for $\sum|Z(A_i)|$ and $\sum |Z'(A_i')|$ for the respective stable
decompositions of all $E\in\kd$. But on each connected component of $\Stab(\kd)$ it is in fact equivalent
to the product metric
$$d(\sigma,\sigma'):=\max\{f(\kp,\kp'),|Z-Z'|\}.$$
As we will restrict to a connected component from the outset, we shall work with this simpler distance function.
Note that due to the definition of $f(\kp,\kp')$, taking into account all objects $E\in\Db(X)$,  the distance between two stability function is difficult to compute explicitly.


\subsection{}\label{subsect:Stand} We shall now specialize to the case that $\kd$ is the bounded derived category $\Db(X):=\Db(\coh(X))$
of the abelian category of coherent sheaves on a complex projective K3 surface. We shall write $\Stab(X)$ for
$\Stab(\Db(X))$.

Stability conditions on higher-dimensional varieties are difficult to construct. On K3 surfaces, Bridgeland constructs 
in \cite{BK3} explicit
examples  of stability conditions as follows. Let $\omega\in\NS(X)_\RR$ be an ample class and let
$B\in\NS(X)_\RR$ be
arbitrary. Consider the linear function $$E\mapsto
Z(E)=\langle\exp(B+i\omega),v(E)\rangle.$$
Here, $v(E)={\rm ch}(E)\sqrt{{\rm td}(X)}\in N(X)\subset H^*(X,\ZZ)$ is the Mukai vector of $E$ and
$\langle~,~\rangle$ is the Mukai pairing.

 Under the additional condition 
that $Z(E)\not\in\RR_{<0}$ for all spherical sheaves (which holds whenever $\omega^2>2$), the function $Z$ has the Harder--Narasimhan property on the  abelian category $\ka(\exp(B+i\omega))$ which is defined as follows (see \cite[Sect.\ 7]{BK3}).

An object $E\in\Db(X)$ is contained in $\ka(\exp(B+i\omega))$ if and only if
there exists an exact triangle
\begin{equation}\label{eqn:disttr}
\kh^{-1}[1]\to E\to\kh^0
\end{equation}
 with
coherent sheaves $\kh^{-1},\kh^0$ satisfying:\\
i) $\kh^{-1}$ is zero or torsion free with $\mu_{\rm max}\leq (B.\omega)$.\\
ii) $\kh^0$ is torsion or $\mu_{\rm min}>(B.\omega)$.

\medskip

The category $\ka(\exp(B+i\omega))$ is the heart of a t-structure that is obtained by tilting the
standard t-structure with respect to the torsion theory described by i) and ii). This defines a stability condition
$\sigma$ depending on $B+i\omega$ whose heart, i.e.\ the abelian category $\kp(0,1]$ of all objects with $\phi^\pm\in(0,1]$,
is precisely $\ka(\exp(B+i\omega))$. We will refer to these stability conditions as \emph{standard stability conditions}.

Standard stability conditions form a  subset $V(X)\subset\Stab(X)$ which via the period map
$\sigma=(\kp,Z)\mapsto Z$ and the Mukai pairing, can be identified with a subset of
$N(X)_\CC$, where $N(X):=H^0\oplus \NS(X)\oplus H^4$ is the algebraic part of $H^*(X,\ZZ)$.
The set of standard stability conditions $V(X)$ can intrinsically be described as follows, see \cite[Prop.\ 10.3]{BK3}.

\begin{prop} {\bf (Bridgeland)}
Suppose $\sigma$ is a stability condition with respect to which for all points $x\in X$ the skyscraper
sheaf $k(x)$ is $\sigma$-stable of phase one. Then $\sigma\in V(X)$.\qqed
\end{prop}

The natural $\widetilde\Gl^+\!\!(2,\RR)$-action on $\Stab(X)$ can be used to describe the set $U(X)$ of all stability conditions with respect
to which all point sheaves $k(x)$ are stable of the same phase. Indeed, $U(X)=V(X)\cdot\widetilde\Gl^+\!\!(2,\RR)$ which can
also be viewed as a $\widetilde\Gl^+\!\!(2,\RR)$-bundle over $V(X)$.

The connected component of $\Stab(X)$ containing $V(X)$ will be denoted $\Sigma$. For $\Sigma$ one has the following
description due to Bridgeland \cite{BK3}.  Consider
the open set $\kp(X)$ of all classes in $N(X)_\CC$ whose real and imaginary part span a positive plane
and let $\kp^+(X)$ be the connected component of $\kp(X)$ that contains
all $\exp(B+i\omega)$ with ample $\omega$.
Then one defines $\kp_0^+(X)$ as the open subset
$\kp^+(X)\setminus\bigcup_{\delta\in\Delta}\delta^\perp$, where $\Delta\subset N(X)$ is the set of $(-2)$-classes.
 
\begin{prop}  {\bf (Bridgeland)}
The period map $\sigma=(\kp,Z)\mapsto Z$ yields a covering map $$\Sigma\to\kp^+_0(X).$$
The group of deck transformation ${\rm Gal}(\Sigma/\kp^+_0(X))$ is naturally identified with the group of
all derived equivalences $\Phi$ preserving $\Sigma$ and acting
trivially on $H^*(X,\ZZ)$.
\end{prop}


\subsection{} For the convenience of the reader we provide the following list of  mostly rather obvious facts on coherent sheaves
on a projective K3 surface $X$. We fix an ample line bundle $\ko(1)$.

\begin{lem}\label{lem:usefulsheaves}
i) If $F$ is a locally free sheaf, then $\Ext^1(\ko(n),F)\cong H^1(X,F^*(n))^*=0$ for $n\gg0$.
ii) If $F\in\coh(X)$ and $\Hom(\ko(n),F)\ne0$ for $n\gg0$, then $F$ contains a non-trivial subsheaf $G\subset F$
with zero-dimensional support.

iii) If $F\in\coh(X)$ is simple, i.e.\ $\End(F)\cong\CC$, then $F$ does not contain a non-trivial proper
subsheaf $0\ne G\subsetneqq F$ with zero-dimensional support.

iv) If $F\in\coh(X)$ is rigid and torsion free, then $F$ is locally free. 
\end{lem}

\begin{proof} Serre duality and Serre vanishing imply i). In order to prove ii), one can argue as follows.
A generic section $t\in H^0(X,\ko(n))$, $n\gg0$, induces naturally an injection
$F(-n)\,\hookrightarrow F$. Thus $h^0(F(-n))\leq h^0(F)$. If indeed $H^0(X,F(-n))=\Hom(\ko(n),F)\ne0$ for $n\gg0$, then we may assume that $h^0(F(-n))=h^0(F)\ne0$ for all $n>0$ (pass to $F(-n_0)$, $n_0\gg0$, if necessary).
Then choose $0\ne s\in H^0(X,F)$ and write it
as $s:\ko_X\twoheadrightarrow\ko_Z\,\hookrightarrow F$ for some non-empty subscheme $Z\subset X$.

For generic $C\in|\ko(n)|$, $n\gg0$ one has exact sequences
$0\to\ko_Z(-n)\to\ko_Z\to\ko_{Z\cap C}\to 0$ and $0\to F(-n)\to F\to F_C\to 0$  with
$\ko_{Z\cap C}\,\hookrightarrow F_C$ and hence $H^0(\ko_{Z\cap C})\subset H^0(F_C)$.
Since $H^0(X,F(-n))=H^0(X,F)$, the section $s\in H^0(\ko_Z)\subset  H^0(F)$ is contained in $H^0(\ko_Z(-n))$. But $H^0(\ko_Z(-n))=0$ for $n>0$ except for $\dim Z=0$.

For iii) consider the non-trivial quotient $F':=F/G$.
If $G$ has zero-dimensional support and $x\in\supp G$, then $k(x)\subset G$. If also  $x\in\supp F'$, then
there exists a surjection $F'\twoheadrightarrow k(x)$ which by composition with $F\twoheadrightarrow F'$
and $k(x)\,\hookrightarrow G\subset F$ yields an endomorphism of $F$ which is not of the form $\lambda\cdot\id$.
If $F'$ and $G$ have disjoint support, then $F\cong F'\oplus G$ which is clearly not simple.

For iv) consider the reflexive hull $F^{**}$ of $F$. The quotient $S$ of $F\subset F^{**}$ is concentrated in dimension zero and
the natural surjection $F^{**}\twoheadrightarrow S$ can be deformed such that $S$ changes its support. Taking kernels yields
a deformation of $F$ which really is non-trivial as the support of its singular part deforms. This contradicts the assumption that
$F$ is rigid. A more explicit dimension count is expressed in \cite[Prop.\ 2.14]{Mu}.
\end{proof}

From these easy facts one can deduce useful information on the heart of a standard stability condition.
Let $\omega\in\NS(X)_\RR$ be an ample class, $B\in\NS(X)_\RR$, and let $\sigma$ be the standard stability condition
with stability function $Z=\langle\exp(B+i\omega,~\rangle$ and heart $\ka:=\ka(B+i\omega)=\kp(0,1]$ (see Section \ref{subsect:Stand}).

\begin{lem}\label{lem:Lemma2}
If $E\in\ka$, then $\Hom(E,\ko(-n)[k])=0$ for $n\gg0$ and $k\leq1$.
\end{lem}

\begin{proof}
By construction, $\ko(-n)\in\ka[-1]$ for $n\gg0$ or, more precisely, for $-n\leq (B.\omega)/(\ko(1).\omega)$.
Hence, $\Hom(\ka,\ko(-n)[k])=0$ for $k\leq0$.
 For $k=1$ use Serre duality to write $$\Hom(E,\ko(-n)[1])=\Ext^1(E,\ko(-n))\cong\Ext^1(\ko(-n),E)^*.$$ Then apply $\Hom(\ko(-n),~)$
 to (\ref{eqn:disttr})  which yields the exact sequence
 $$\Ext^2(\ko(-n),\kh^{-1})\to\Ext^1(\ko(-n),E)\to\Ext^1(\ko(-n),\kh^0).$$
 Then for $n\gg0$ Serre vanishing  yields $\Ext^2(\ko(-n),\kh^{-1})=H^2(X,\kh^{-1}(n))=0$ and similarly $ \Ext^1(\ko(-n),\kh^0)=H^1(X,\kh^0(n))=0$.
 \end{proof}

A similar `dual' statement for spherical objects will be proved in Lemma \ref{lem:Lemma1}.
\section{Spherical objects}\label{sect:Spherical}


\subsection{} Let us recall the definition of a spherical object. We shall work with a \emph{K3 category} $\kd$
which later will be $\Db(X)$, the bounded derived category of coherent sheaves on
a K3 surface $X$.  Recall that a K3 category is a linear triangulated category of finite type with the shift
$[2]$ defining a Serre functor.

\begin{definition} An object $A\in\kd$ is called \emph{spherical} if
$$\Ext^*(A,A)\cong H^*(S^2,\CC).$$ 
By $\ks\subset{\rm Ob}( \kd)$ we denote the collection of all spherical objects in $\kd$.
\end{definition}


Thus  $A\in \ks$ if and only if $A$ is simple (i.e.\ $\End(A)\cong\CC$), rigid (i.e.\ $\Ext^1(A,A)=0$),
and  $\Ext^i(A,A)=0$ for $i<0$. The easiest examples 
are provided by line bundles on a K3 surface $X$ viewed as objects in the K3 category $\Db(X)$.

\smallskip

To shorten the notation we will sometimes write $\ext^i(A,B)=\dim\Ext^i(A,B)$. The following results go back to Mukai, see e.g.\
\cite[Cor.\ 2.8]{Mu}. In this form they can be found in  \cite[Lem.\ 2.7, Prop.\ 2.9]{HMS} (see also \cite[Lem.\ 12.2]{BK3}).

\begin{lem}\label{lem:crucial}
 Consider in the K3 category $\kd$ an exact triangle
$$\xymatrix{
A\ar[r]^i&E\ar[r]^j&B\ar[r]^\delta&A[1]}$$ such that
\[
\Ext^r(A,B)=\Ext^s(B,B)=0~\text{for~} r\leq0\text{~and~} s<0.
\]
Then
\[
\ext^1(A,A)+\ext^1(B,B)\leq\ext^1(E,E).
\]
\end{lem}

The following two consequences hold true for arbitrary slicings, no stability function is needed.

\begin{cor}\label{cor:HMS1} Let $\sigma$ be a stability condition on $\kd$ and $A\in\ks$. If $A_1,\ldots,A_k$ are the
 $\sigma$-stable factors  of $A$  (cf.\ Section \ref{sect:basics}), then $A_1,\ldots,A_k\in\ks$.\qqed
\end{cor}

An object $E\in\kd$ is called \emph{semirigid} if $\Ext^1(E,E)$ is two-dimensional. If $x\in X$ is a closed point of
a K3 surface $X$, then $k(x)$ is a semirigid object in $\Db(X)$.

\begin{cor}\label{cor:HMS2}
Let $\sigma$ be a stability condition on a K3 category $\kd$ and let $E$ be a semirigid object. Then
the $\sigma$-stable factors $E_1,\ldots, E_k$ of $E$ are spherical or semirigid. In fact, at most one $E_i$ can be semirigid.
\qqed
\end{cor}

\subsection{}

Consider two stability conditions $\sigma$ and $\sigma'$ on the K3 category $\kd$. 
The proof of the following result only uses the underlying slicings, $\kp$ resp.\ $\kp'$,
and the property that all $\kp(\phi)$ and $\kp'(\phi)$ are abelian.

\begin{prop}\label{prop:equivcond}
The following conditions are equivalent:

i) For all $\phi\in\RR$ one has $\kp(\phi)^s\cap\ks=\kp'(\phi)^s\cap\ks$.

ii)  For all $\phi\in\RR$ one has $\kp(\phi)\cap\ks=\kp'(\phi)\cap\ks$.

iii) For all $A\in\ks$ one has $\phi^\pm(A)={\phi'}^\pm(A)$.
\end{prop}

\begin{proof}
Assume iii). An object $E$ is $\sigma$-semistable of phase $\phi$ if and only of $\phi^+(E)=
\phi^-(E)=\phi$. But for $A\in\ks$ this is, assuming iii), equivalent to ${\phi'}^+(A)=
{\phi'}^-(A)=\phi$.  Hence such an $A$ is also $\sigma'$-semistable of the same phase $\phi$. Thus, ii) holds.

Assume ii). If $A\in\ks$ is $\sigma$-stable of phase $\phi$, then $A$ is in particular $\sigma$-semistable
of phase $\phi$ and hence by ii) also $\sigma'$-semistable of phase $\phi$. If $A$ is not $\sigma'$-stable, then there
exists a minimal proper subobject $A'\subset A$ in the abelian category $\kp'(\phi)$. Then $A'$ is $\sigma'$-stable and as a stable factor of a spherical object,
$A'$ is also spherical  (cf.\ Corollary \ref{cor:HMS1}). Hence by ii), $A'$ is also $\sigma$-semistable of phase $\phi$.
One would now like to argue that then the inclusion $A'\subset A$ in $\kp'(\phi)$ must be an isomorphism because $A$ was $\sigma$-stable.
However a priori we do not know that $A'\to A$ is still an injection in $\kp(\phi)$. But
since $A\in\kp(\phi)$ is $\sigma$-stable, $A'\in\kp(\phi)$,  and $A'\to A$ is non-trivial, $A$ is a $\sigma$-stable factor of $A'$ and
$A'\to A$ is a surjection in $\kp(\phi)$. The  $\sigma$-stable factors of its kernel in $\kp(\phi)$ are
$\sigma$-stable factors of the spherical $A'$ and hence also spherical. 
Thus the short exact sequence $0\to {\rm Ker}\to A'\to A\to 0$ in $\kp(\phi)$ is also a short exact sequence in $\kp'(\phi)$, but as $A'$ was a subobject
of $A$ in $\kp'(\phi)$ this shows
${\rm Ker}=0$. Hence $A'\cong A$ and thus $A$ is $\sigma'$-stable of phase $\phi$.
This shows i).

Assume i). Consider a $\sigma$-stable filtration $F_1\to\ldots\to F_n=A$ with $\sigma$-stable factors $A_i$ of
phase $\phi_i$. Since $A\in\ks$,  all $A_i\in \ks$. Hence, all $A_i$ are by i) also $\sigma'$-stable of phase $\phi_i$.
In particular, the given filtration is also a  stable filtration with respect to $\sigma'$. But then $\phi^+(A)=\phi_1
={\phi'}^+(A)$ and $\phi^-(A)=\phi_n
={\phi'}^-(A)$. This shows iii).
\end{proof}


\subsection{}
In analogy to Lemma \ref{lem:Lemma2} one has the following result for spherical objects
in the bounded derived category $\Db(X)$ of a complex projective K3 surface $X$. As before,
$\ka$ is the heart of a standard stability condition
with stability function $Z=\langle\exp(B+i\omega),~\rangle$.

\begin{lem}\label{lem:Lemma1}
If $A\in\ka$ is spherical, then $\Hom(\ko(n),A[k])=0$ for all $k\leq0$ and $n\gg0$.
\end{lem}

\begin{proof}
By stability, $\Hom(\ka,\ka[k])=0$ for $k<0$. Since $\ko(n)\in\ka$ for $n\gg0$, or more precisely for $n>(B.\omega)/(\ko(1).\omega)$, this proves the vanishing for
negative $k$. To prove the vanishing for $k=0$ apply $\Hom(\ko(n),~)$ to (\ref{eqn:disttr}) for $A$ which yields the exact sequence
$$\Ext^1(\ko(n),\kh^{-1})\to\Hom(\ko(n),A)\to\Hom(\ko(n),\kh^0).$$
As $A$ is spherical and $\Hom(\kh^{-1}[1],\kh^0)=0$, Lemma \ref{lem:crucial} shows that $\kh^{-1}$ and $\kh^0$ are both rigid sheaves.
Thus $\kh^{-1}$ is a rigid torsion free sheaf and therefore locally free (see Lemma \ref{lem:usefulsheaves}, iv)). By
Lemma \ref{lem:usefulsheaves}, i) one finds $\Ext^1(\ko(n),\kh^{-1})=0$ for $n\gg0$. Thus, if $\Hom(\ko(n),A)\ne0$ for
$n\gg0$, then $\Hom(\ko(n),\kh^0)\ne0$ for $n\gg0$. By Lemma \ref{lem:usefulsheaves}, ii), this means that the zero-dimensional part
 $G:=T_0(\kh^0)\subset \kh^0$ of $\kh^0$ is non-trivial. If $\kh^0$ is not only rigid but in fact spherical, then Lemma \ref{lem:usefulsheaves}, iii) would
show that $\kh^0$ is zero-dimensional and in fact $\kh^0\cong k(x)$. Clearly,  the latter
would contradict rigidity of $\kh^0$. If $\kh^0$ is rigid but not simple, one can argue as follows.
Note that $\Ext^i(G,\kh^0/G)=0$ for $i\leq0$ and $\Ext^i(\kh^0/G,\kh^0/G)=0$ for $i<0$. Then by Lemma  \ref{lem:crucial} one finds
$\ext^1(G,G)+\ext^1(\kh^0/G,\kh^0/G)\leq \ext^1(\kh^0,\kh^0)=0$, but clearly the zero-dimensional sheaf $G$ deforms and hence
$\Ext^1(G,G)\ne0$ which yields a contradiction.
\end{proof}


\section{Stability conditions via spherical objects}\label{sect:Main}

Let $X$ be a complex projective K3 surface and $\Sigma\subset\Stab(X)$  the distinguished
connected component of the space of locally finite numerical stability conditions
on $\Db(X)$ (see \cite{BK3}).

This section is entirely devoted to the proof of the following

\begin{thm}\label{thm:main}
Suppose $\sigma=(\kp,Z),\sigma'=(\kp',Z')$ are stability conditions in
$\Sigma$. Then $$\sigma=\sigma'$$ if and only if 
$Z=Z'$ and for every spherical object
$A\in\Db(X)$:
\begin{equation}\label{eqn:sphsst}
A \textrm{ is } \sigma\textrm{-semistable of phase } \varphi \textrm{ if and only if } A \textrm{ is } \sigma'\textrm{-semistable of phase }\varphi.
\end{equation}
\end{thm}

As we shall see, the proof really uses that both stability conditions, $\sigma$ and $\sigma'$, are contained
in the distinguished component $\Sigma$ or, slightly weaker,  that one of the two is contained in $\Sigma$
and that the set of all point sheaves $k(x)$ is of bounded mass with respect to the other.

\smallskip

The proof proceeds in three steps. We shall first assume  that $\sigma$ is a standard stability condition (Section \ref{subsect:pf1})
and then reduce to this case by applying autoequivalences. The case that $\sigma$ can only be transformed into a stability condition that
is a limit of standard stability conditions will be dealt with in Section \ref{subsect:pf3}

We will frequently use the observation (see Proposition \ref{prop:equivcond}) that  (\ref{eqn:sphsst})
is equivalent to:
\begin{equation}\label{eqn:sphst}
A \textrm{ is } \sigma\textrm{-stable of phase } \varphi \textrm{ if and only if } A \textrm{ is } \sigma'\textrm{-stable of phase }\varphi.
\end{equation}

\subsection{}\label{subsect:pf1}
Assume $\sigma$ and $\sigma'$ satisfy $Z=Z'$ and (\ref{eqn:sphsst}) (or, equivalently, (\ref{eqn:sphst})) and that
in addition $\sigma\in V(X)$. In particular 
$Z=\langle\exp(B+i\omega),~\rangle$ for some ample $\omega$ 
and all point sheaves $k(x)$ are $\sigma$-stable of phase one. In order to show that $\sigma=\sigma'$ 
with  $\sigma'$ as in Theorem \ref{thm:main}, it suffices
to show that the point sheaves $k(x)$ are also $\sigma'$-stable of phase one (cf.\ Section \ref{subsect:Stand}).
For this, the assumption that $\sigma'$ is contained
in the connected component $\Sigma$ is not needed.

To shorten the notation we shall denote the heart $\ka(B+i\omega)=\kp(0,1]$ of $\sigma$ simply by $\ka$.

\begin{lem}\label{lem:phaseone}
If $k(x)$ is $\sigma'$-stable, then its phase with respect to $\sigma'$ is one, i.e.\ $\phi'(k(x))=1$.
\end{lem}

\begin{proof}
Pick a line bundle $L$ with $(L.\omega)>(B.\omega)$. Then $L\in\ka$ by definition of $\ka=\ka(B+i\omega)$.
The line bundle $L$ is a spherical object and by Corollary \ref{cor:HMS1} all $\sigma$-stable factors $L_i$ of $L$ are spherical as well.
Since $L\in\ka$, their phases satisfy $\phi(L_i)\in(0,1]$. 

By our assumption  on $\sigma'$ (see (\ref{eqn:sphst})), the $L_i$ are then also $\sigma'$-stable of phase $\phi'(L_i)=\phi(L_i)\in(0,1]$.
 Clearly, any line bundle $L$ admits a non-trivial morphism $L\to k(x)$ and
hence at least one of the $\sigma$-stable factors $L_i$ admits a non-trivial morphism
$L_i\to k(x)$. Since we assume $k(x)$ to be $\sigma'$-stable, its $\sigma'$-phase is well defined
and thus satisfies $0<\phi'(L_i)\leq\phi'(k(x))$. On the other hand, by Serre duality,
$\Hom(L_i, k(x))\ne0$ implies $\Hom(k(x), L_i[2])\ne0$. The latter yields
$\phi'(k(x))\leq\phi'(L_i[2])\leq 3$. Moreover, $\phi'(k(x))=\phi'(L_i)=3$ can only occur if the two $\sigma'$-stable objects
$k(x)$ and $L_i[2]$ are isomorphic, which is absurd as one is semirigid and the other is spherical. Thus, $\phi'(k(x))\in(0,3)$. As
$Z=Z'$ and $Z(k(x))=-1$, this readily shows $\phi'(k(x))=1$.
\end{proof}

Suppose $k(x)$ is not $\sigma'$-stable. Then there exists a $\sigma'$-stable decomposition, i.e.\ a diagram
$$\xymatrix{&F_1\ar@{=}[dl]\ar[rr]&&F_2\ar[dl]\ar[r]&\ldots\ar[r]& F_{m-1}\ar[rr]&&F_m\ar@{=}[r]\ar[dl]&k(x)\\
A_1&&A_2\ar[ul]^{[1]}&&&&A_m\ar[ul]^{[1]}&&&}$$
where the $A_i$ are $\sigma'$-stable with $\phi'(A_1)\geq\ldots \geq \phi'(A_m)$ and $m>1$.
By Corollary \ref{cor:HMS2} at most one $A_i$ is not spherical and if there is one, it is semirigid.

\smallskip

{\bf i)} Suppose $A_1$ and $A_m$ are both  spherical. Then by (\ref{eqn:sphst}), both are also $\sigma$-stable and for their phases one has
$\phi(A_1)=\phi'(A_1)$ and $\phi(A_m)=\phi'(A_m)$. Since $A_1=F_1\to k(x)$ is not trivial, $\phi(A_1)\leq\phi(k(x))=1$ and equality would
imply $A_1=k(x)$ which can be excluded as in the proof of Lemma \ref{lem:phaseone}.
Similarly, $k(x)=F_m\to A_m$ is not trivial and hence $1=\phi(k(x))\leq\phi(A_m)$.
This yields the contradiction $1>\phi(A_1)=\phi'(A_1)\geq\phi'(A_m)=\phi(A_m)\geq1$.

\smallskip

{\bf ii)} Suppose $A_m$ is semirigid. Then $A_1,\ldots,A_{m-1}$ are spherical and hence also $\sigma$-stable with phases
$\phi(A_i)=\phi'(A_i)$, $i=1,\ldots,m-1$. As above, $\Hom(A_1,k(x))\ne0$ implies
$1>\phi(A_1)=\phi'(A_1)\geq\ldots\geq\phi(A_{m-1})=\phi'(A_{m-1})$. Thus $A_i\in\ka[k_i]$, $i=1,\ldots, m-1$, with $k_i\leq0$.
Then Lemma \ref{lem:Lemma1} shows $\Hom(\ko(n),A_i)=0$ for $i=1,\ldots,m-1$ and $n\gg0$.
Since $\Hom(\ko(n),k(x))\ne0$ for all $n$, we find
$\Hom(\ko(n),A_m)\ne0$ for $n\gg0$.

 Clearly, $\ko(n)\in\ka$ for $n\gg0$ and therefore all $\sigma$-stable factors
$L_i$ of $\ko(n)$, which by Corollary \ref{cor:HMS1} are also spherical, have phases $\phi(L_i)\in(0,1]$.
By (\ref{eqn:sphst}) the $L_i$ are also $\sigma'$-stable with $\sigma'$-phases 
$\phi'(L_i)=\phi(L_i)$. Since $\Hom(\ko(n),A_m)\ne0$ implies $\Hom(L_i,A_m)\ne0$ for at least one $L_i$, stability yields
$0<\phi(L_i)=\phi'(L_i)\leq\phi'(A_m)$.

Thus one finds $1>\phi'(A_1)\geq\ldots\geq\phi'(A_{m-1})\geq \phi'(A_m)\geq0$. Since $Z'(k(x))=\sum Z'(A_i)$ and $Z'(k(x))=Z(k(x))=-1$,
this is impossible.

\smallskip

{\bf iii)} Suppose $A_1$ is semirigid. Then $A_2,\ldots,A_{m}$ are spherical and hence also $\sigma$-stable with phases
$\phi(A_i)=\phi'(A_i)$, $i=2,\ldots,m$. Using that $k(x)\to A_m$ is non-trivial and not an isomorphism, one finds
$1<\phi(A_m)=\phi'(A_m)\leq\ldots\leq\phi'(A_2)=\phi(A_2)$. Hence, $A_i\in\ka[k_i]$, $i=2,\ldots,m$, with $k_i\geq1$.
Then Lemma \ref{lem:Lemma2} shows $\Hom(A_i,\ko(-n)[2])=0$ for $i=2,\ldots,m$ and $n\gg0$. Since
$\Hom(k(x),\ko(-n)[2])=\Hom(\ko(-n),k(x))^*\ne0$ for all $n$, this yields $\Hom(A_1,\ko(-n)[2])\ne0$ for $n\gg0$.

For a fixed such $n$, consider the $\sigma$-stable factors $L_i$ of $\ko(-n)[2]$ which are contained in $\ka[1]$ and hence
$\phi(L_i)\in(1,2]$. Again, the $L_i$ are spherical (cf.\ Corollary \ref{cor:HMS1}) and hence by (\ref{eqn:sphst}) also $\sigma'$-stable of phase $\phi'(L_i)=\phi(L_i)$. Since
$\Hom(A_1,\ko(-n)[2])\ne0$ implies  $\Hom(A_1,L_i)\ne0$ for at least one $L_i$, stability yields
$\phi'(A_1)\leq 2$

Thus one finds $1<\phi'(A_m)\leq\ldots\leq\phi'(A_1)\leq2$. As above, this contradicts $Z'(k(x))=Z(k(x))=-1$.

This concludes the proof of Theorem \ref{thm:main} in the case that $\sigma\in V(X)$. In Section \ref{subsect:pf3}
we will use similar arguments for the case that $\sigma\in\partial V(X)$, but they will have to be applied
to small deformations of $\sigma$ and $\sigma'$ which makes it more technical.


\subsection{}
Suppose now that $\sigma,\sigma'\in\Sigma$ satisfy $Z=Z'$ and (\ref{eqn:sphsst}) (or, equivalently, (\ref{eqn:sphst})).
In order to show that then $\sigma=\sigma'$ it suffices
to find an autoequivalence $\Phi\in\Aut(\Db(X))$ such that $\Phi(\sigma)=\Phi(\sigma')$. Since the set of spherical objects $\ks\subset{\rm Ob}(\Db(X))$
is invariant under the action of $\Aut(\Db(X))$, the new stability conditions $\Phi(\sigma),\Phi(\sigma')$ still satisfy (\ref{eqn:sphsst}).  

Recall that for any $\sigma\in\Sigma$ there exists  $\Phi\in\Aut(\Db(X))$, such that $\Phi(\sigma)$ is contained in the closure $\overline{U(X)}$ of the
open set $U(X)\subset\Sigma$ of all stability conditions with respect to which all point sheaves $k(x)$ are stable of the same phase (see \cite{BK3}).
Moreover, $U(X)$ is a principal $\widetilde\Gl^+\!\!(2,\RR)$-bundle over $V(X)\subset U(X)$ (see \cite[Sect.\ 11]{BK3} or Section \ref{subsect:Stand}).

Thus, if $\Phi$ can be found such that $\Phi(\sigma)\in U(X)$ (and not only in its closure), then there exists a $g\in\widetilde \Gl^+\!\!(2,\RR)$ with
$g^{-1}\Phi(\sigma)\in V(X)$. By (\ref{subsect:pf1}), applied to $g^{-1}\Phi(\sigma)$ and $g^{-1}\Phi(\sigma')$, one concludes $g^{-1}\Phi(\sigma)=g^{-1}\Phi(\sigma')$ and hence $\sigma=\sigma'$.


\subsection{}\label{subsect:pf3}
Eventually we have to deal with the case that one only finds  $\Phi\in\Aut(\Db(X))$ such that
$\Phi(\sigma)$ is in the boundary of $U(X)$. By applying an appropriate $g\in\widetilde\Gl^+\!\!(2,\RR)$ we can reduce to the
case that  $\sigma\in\partial V(X)$, i.e.\ all $k(x)$ are $\sigma$-semistable (but not necessarily stable) of phase one,  and
$\sigma'\in \Sigma$. 

Pick a path $\sigma_t$, $0\leq t\ll1$ with $\sigma_0=\sigma$ and $\sigma_t\in V(X)$ for $t>0$. The stability function of $\sigma_t$ shall
be denoted  $Z_t$ and  for a $\sigma_t$-semistable object $B$ its phase is  $\phi_t(B)$.
Since $Z=Z'$ and $\sigma'\in\Sigma$, the path $\sigma_t$ (or rather its image in $\kp^+_0(X)$) can be lifted
uniquely to a path $\sigma'_t$ in $\Sigma$ with $\sigma'_0=\sigma'$. Then, by construction, the stability function $Z_t'$ of $\sigma_t'$ equals $Z_t$.
The phase of a $\sigma_t'$-semistable object $B$ shall be denoted $\phi_t'(B)$. 

In the following, $\sigma_t$-semistability  of an object will
mean semistability  for all small $t$ (depending on the object) and similarly for $\sigma'_t$-semistability.
Note that semistability is a closed condition, so semistability
for all small $t>0$ will imply semistability for $t=0$. The same does not hold for semistability replaced by
stability. So, when we say an object is $\sigma_t$-stable, it means that it is $\sigma_t$-stable for all small
$t>0$. The latter implies that it is also $\sigma$-semistable, but not necessarily $\sigma$-stable.

We continue to assume (\ref{eqn:sphsst}) (or, equivalently, (\ref{eqn:sphst})) for the two stability conditions
$\sigma$ and $\sigma'$.  The condition is preserved under small deformation as shown by the following

\begin{lem}\label{prop-staysst} Suppose $A$ is a spherical object. Then the path $\sigma_t$ can be chosen such that
$A$ is $\sigma_t$-semistable if and only if $A$ is $\sigma_t'$-semistable. Moreover, in this case $\phi_t(A)=\phi_t'(A)$.
\end{lem}

\begin{proof} Recall  that for fixed $Z\in\kp^+(X)$ and an
arbitrary norm on $N(X)_\RR$ there exists a
constant $C$ such that for all $(-2)$-classes $\delta\in N(X)$ one has $\|\delta\|^2\leq 2(1+C|Z(\delta)|^2)$.
This can be found implicitly in the proof of \cite[Lem.\ 8.1]{BK3} (and explicitly in the first version of the paper).
Hence the set of $(-2)$-classes $\delta\in N(X)$ with bounded $Z(\delta)$ is finite. 

Therefore it suffices to prove that the assertion holds for $A$ once it holds for all spherical objects $B$ with $|Z(B)|<|Z(A)|$.
A priori the interval $t\in[0,\varepsilon)$ for which semistability with respect to $\sigma_t$ resp.\ $\sigma'_t$ coincide
can get smaller when passing from $A$ to $B$. But only finitely many steps are necessary and, as we shall
see,  at each step only finitely many spherical objects are  involved. 

Suppose $A$ is $\sigma_t$-semistable but not $\sigma'_t$-semistable for $t>0$. Then there exists
a $\sigma'_t$-stable decomposition of $A$ with $\sigma_t'$-stable factors $B_i$ such that $\phi_t'(B_1)\geq\ldots\geq\phi_t'(B_k)$.
The arguments to show this can be found in the proof of \cite[Prop.\ 9.3]{BK3}. Take a compact neighbourhood $K$ of $\sigma'$ and consider the
set $T(A,K)$ of all objects $B$ that are stable factors of $A$ with respect to some $\sigma_t'\in K$.  This set is of bounded mass
and by \cite[Prop. 9.3]{BK3} there is a finite chamber structure of $K$ such that
for an object $B\in T(A,K)$ (semi)stability is constant within a chamber. This chamber structure can be refined such that within one chamber 
$\log(\phi_t'(B_1)/\phi'_t(B_2))$ does not change signs for all $B_1,B_2\in T(A,K)$. By the finiteness of the set of Mukai vectors $\{v(B)~|~B\in T(A,K)\}$
(see  \cite[Lem.\ 9.2]{BK3}) the new chamber structure is still finite. Hence $\sigma'$ will be in the closure of one chamber and we choose
$\sigma'_t$ in this chamber and find the stable decomposition as claimed.

Note that by the assumption that $A$ is not $\sigma_t'$-semistable, one has $\phi_t'(B_1)>\phi_t'(B_k)$ for $t>0$.
Since $A$ is spherical, also its stable factors $B_i$ are spherical (cf.\ Corollary \ref{cor:HMS1}). For $t=0$ one has
$\phi(B_1)=\phi'(B_1)=\ldots=\phi'(B_k)=\phi(B_k)$, because $A$ is $\sigma'$-semistable by (\ref{eqn:sphsst}).
Hence $Z'(B_i)=Z(B_i)\in Z(A)\RR_{>0}$. Since $Z(A)=\sum Z(B_i)$, one has
$|Z(B_i)|< |Z(A)|$. But then the assertion of the lemma  holds for the $B_i$ which are $\sigma_t'$-semistable. 
(At this point the path $\sigma_t$ has to be adjusted to work for the $B_i$ as well. As mentioned earlier, this procedure
really works, because only finitely many objects are eventually used.)
Thus, the $B_i$ are $\sigma_t$-semistable with $\phi_t(B_i)=\phi_t'(B_i)$. Hence
$\phi_t(B_1)=\phi_t'(B_1)>\phi_t'(B_k)=\phi_t(B_k)$ for $t>0$ contradicting the $\sigma_t$-semistability of $A$.

If $A$ is semistable with respect to $\sigma_t$ and $\sigma'_t$, then $\phi(A)=\phi'(A)$ by (\ref{eqn:sphsst}). As $Z_t=Z'_t$, this yields
 $\phi_t(A)=\phi_t'(A)$.
\end{proof}

Let us now turn to the proof of Theorem \ref{thm:main} in this situation. Morally, Lemma  \ref{prop-staysst} says 
that we can apply Section \ref{subsect:pf1} to the stability conditions $\sigma_t$ and $\sigma'_t$ for some small $t>0$.
However, the chamber structure that takes care of all the objects involved may not be locally finite near $\sigma'$. Indeed,
one would start with the $\sigma'_t$-stable factors $A_i$ of some $k(x)$ and in the next step would need to consider
the $\sigma_t$-stable
factors of the $A_i$ and so forth. So we have to run the arguments of Section  \ref{subsect:pf1} once more while keeping track of the
deformation to the interior of $V(X)$ (which makes everything  more technical).

We shall prove that each $k(x)$ is $\sigma'_t$-stable
of phase one for small $t>0$. Since the family of all $k(x)$ is of bounded mass in $\Sigma$, this suffices to conclude
that $\sigma'\in\partial V(X)$. Indeed by \cite[Prop.\ 9.3]{BK3} the chamber structure of a compact neighbourhood of $\sigma'$ with
respect to $\{k(x)\}$ is finite and hence  all $k(x)$ will be $\sigma_t'$-semistable for $t$ small but independent of the particular point sheaf $k(x)$.
Moreover, as in Section \ref{subsect:pf1}, the phase will be one and hence $\sigma_t'$ is a standard stability condition. Then $Z_t'=Z_t$ and the fact that a standard
stability condition is determined by its stability function shows $\sigma_t=\sigma_t'$ and hence $\sigma=\sigma'$.

\medskip

Suppose $k(x)$ is not $\sigma_t'$-stable.
Then there exists a  decomposition as in Section \ref{subsect:pf1} with factors $A_1,\ldots, A_m$ which are
$\sigma_t'$-stable and satisfy $\phi'_t(A_1)\geq\ldots\geq\phi'_t(A_m)$. This follows from \cite[Prop.\ 9.3]{BK3} (see also the arguments in the proof of Lemma \ref{prop-staysst}).
Note that then $A_1,\ldots,A_m$ are still $\sigma'$-semistable but not necessarily $\sigma'$-stable.
In the following, we use similar  arguments as in Section \ref{subsect:pf1}. In particular, we  distinguish three cases.

\smallskip

{\bf i)} Suppose $A_1$ and $A_m$ are both spherical. Then by (\ref{eqn:sphsst}) they are also $\sigma$-semistable
with $\phi(A_1)=\phi'(A_1)$ and $\phi(A_m)=\phi'(A_m)$. Due to the existence of the non-trivial morphisms $A_1\to k(x)$ and
$k(x)\to A_m$ and the $\sigma$-semistability of $k(x)$, this yields
$\phi'(A_1)=\phi(A_1)\leq \phi(k(x))=1$ and $1=\phi(k(x))\leq\phi(A_m)=\phi'(A_m)$. Together with $\phi'_t(A_1)\geq\ldots\geq\phi'_t(A_m)$ 
one finds that $k(x)$ is $\sigma'$-semistable. In fact more is true. Since $k(x)$ is $\sigma_t$-stable for $t>0$ and by Lemma
\ref{prop-staysst} $A_1$ and $A_m$ are $\sigma_t$-semistable with $\phi_t(A_i)=\phi_t'(A_i)$, one obtains
$1\geq\phi_t(A_1)\geq\phi_t(A_m)\geq1$. As the $\sigma_t$-stable semirigid $k(x)$ cannot  be a $\sigma_t$-stable factor of the
spherical $A_1$ (use Corollary \ref{cor:HMS1}),
the first inequality must be strict which is absurd for $m>1$. Thus, $k(x)$ is $\sigma_t'$-stable for $t>0$ of phase one. Hence, if 
we are in case {\bf i)} for all $x\in X$, then $\sigma_t\in\partial V(X)$.

\smallskip

{\bf ii)} Suppose $A_m$ is semirigid. Then $A_1,\ldots,A_{m-1}$ are spherical and by  (\ref{eqn:sphsst})
also $\sigma$-semistable of phase $\phi(A_i)=\phi'(A_i)$. The existence of the non-trivial $A_1\to k(x)$ and 
$\sigma$-semistability of $k(x)$ yield
$1\geq\phi(A_1)=\phi'(A_1)\geq\ldots\geq\phi'(A_{m-1})=\phi(A_{m-1})$.

By Lemma \ref{prop-staysst} the $A_i$, $i=1,\ldots,m-1$, are $\sigma_t$-semistable of phase
$\phi_t(A_i)=\phi_t'(A_i)$. Thus, $\phi_t(A_1)=\phi_t'(A_1)\geq\ldots\geq\phi_t'(A_{m-1})=\phi_t(A_{m-1})$. Moreover,
$\sigma_{t>0}$-stability of $k(x)$ implies $1=\phi_t(k(x))\geq\phi_t(A_1)$ for $t>0$. (Actually,  $\phi_t(k(x))=\phi_t(A_1)$ can
be excluded for $t>0$, because as above  the semirigid $k(x)$ cannot be a stable factor of the spherical $A_1$, see Corollary
\ref{cor:HMS1}). Thus $A_1,\ldots,A_{m-1}$ are $\sigma_{t>0}$-semistable of phase
$\leq1$ (in fact, $<1$) and by Lemma \ref{lem:Lemma1} this proves $\Hom(\ko(n),A_i)=0$ for $n\gg0$ and $i=1,\ldots,m-1$. As in
Section \ref{subsect:pf1}, ii) this yields $\Hom(\ko(n),A_m)\ne0$ for $n\gg0$. 

Let now $L_1,\ldots,L_k$ be the $\sigma$-stable factors of $\ko(n)$ with $\phi(L_1)\geq\ldots\geq\phi(L_k)$.
They are again spherical and hence by (\ref{eqn:sphsst})
also $\sigma'$-stable. Stability is an open property for objects with primitive Mukai vector (see \cite[Prop.\ 9.4]{BK3}). 
Hence the $L_i$ are stable with respect to $\sigma_t$ and $\sigma_t'$ and, moreover,
$\phi_t(L_i)=\phi_t'(L_i)$. For $t>0$ all stable factors of $\ko(n)$, $n\gg0$, have phases
in $(0,1]$. The non-trivial $\ko(n)\to L_k$ and $\sigma_t$-stability of $L_k$
therefore imply $\phi_t(L_k)> 0$ for $t>0$.
In the limit we still have $\phi(L_k)\geq0$ and hence $\phi(L_i)\geq0$ for all $L_i$.

Since $\Hom(L_i,A_m)\ne0$ for at least one $L_i$, semistability of $A_m$ and $L_i$ with respect to $\sigma'$ yields $\phi'(L_i)\leq\phi'(A_m)$ and hence $0\leq\phi'(A_m)\leq\ldots\leq\phi'(A_1)\leq 1$. (The last inequality is
a priori not strict.) This  contradicts $Z'(k(x))=-1$ except for the case that $\phi'(A_m)=\ldots=\phi'(A_1)=1$.
However, if $\phi'(A_m)=1$, then for small $t>0$ one still has $\phi'_t(A_m)>0$ and
thus $0<\phi'_t(A_m)\leq\ldots\leq\phi'_t(A_1)< 1$ where the last inequality is strict for $t>0$.
This contradicts $Z'_t(k(x))=-1$.

\smallskip

{\bf iii)} Suppose $A_1$ is semirigid. Then $A_2,\ldots,A_m$ are spherical and hence by (\ref{eqn:sphsst}) also $\sigma$-semistable
of phase $\phi(A_i)=\phi'(A_i)$. The existence of the non-trivial $k(x)\to A_m$ implies
$\phi(A_2)=\phi'(A_2)\geq\ldots\geq\phi'(A_m)=\phi(A_m)\geq\phi(k(x))=1$. By Lemma \ref{prop-staysst}
we know that $A_2,\ldots,A_m$ are $\sigma_t$-semistable of phase $\phi_t(A_i)=\phi'_t(A_i)$ and 
thus $\phi_t(A_2)=\phi'_t(A_2)\geq\ldots\geq\phi'_t(A_m)=\phi_t(A_m)>\phi(k(x))=1$. The 
last inequality is strict because the   $\sigma_{t>0}$-stable semirigid $k(x)$ cannot be a stable factor of the
spherical $A_m$ (Corollary \ref{cor:HMS1}).

By Lemma \ref{lem:Lemma2} one then has $\Hom(A_i,\ko(-n)[2])=0$ for $n\gg0$ and $i=2,\ldots,m$. And as in Section \ref{subsect:pf1}, iii) this yields $\Hom(A_1,\ko(-n)[2])\ne0$. Consider the $\sigma$-stable factors $L_1,\ldots,L_k$
of $\ko(-n)[2]$ with $\phi(L_1)\geq\ldots\geq\phi(L_k)$. Since they are spherical, they are also $\sigma'$-stable and 
hence $\sigma'_t$-stable, as stability is open for objects with primitive Mukai vector by \cite[Prop.\ 9.4]{BK3}.
Using
Lemma \ref{prop-staysst} one finds that the $L_i$ are semistable with respect to $\sigma_t'$ and $\sigma_t$. Moreover,
$\phi_t'(L_i)=\phi_t(L_i)$. Using that all $\sigma_t$-stable factors of $\ko(-n)[2]$ have phases in $(1,2]$
and the existence of the non-trivial $L_1\to \ko(-n)[2]$, one finds $\phi_t(L_1)\leq2$.
Thus also $\phi(L_1)\leq2$ and hence $\phi(L_k)\leq\ldots\leq\phi(L_1)\leq2$.

As $\Hom(A_1,L_i)\ne0$ for at least one $L_i$ and both, $A_1$ and $L_i$, are $\sigma'$-stable, one finds
$1\leq\phi'(A_m)\leq\ldots\leq\phi'(A_1)\leq2$. The latter contradicts $Z'(k(x))=-1$ except for the case
that $\phi'(A_m)=\ldots=\phi'(A_1)=1$. However, if $\phi'(A_1)=1$, then for small $t>0$ still $\phi'_t(A_1)<2$
and hence $1<\phi'_t(A_m)\leq\ldots\leq\phi'_t(A_1)<2$ where the first inequality is strict for $t>0$.
This, once more, contradicts $Z_t'(k(x))=-1$.

\begin{remark}
The rough idea of the above arguments goes as follows.
 If for two stability conditions $\sigma$ and $\sigma'$ with the same stability function
$Z=Z'$ an object $E$ is stable with respect to $\sigma$ but not with respect to $\sigma'$, then pass to the $\sigma'$-stable factors $A_i$ of $E$.
Either for all $A_i$ one has $\ext^1(A_i,A_i)<\ext^1(E,E)$ or all $A_i$ are spherical except for one, say $A_{i_0}$, for
which $\ext^1(A_{i_0},A_{i_0})=\ext^1(E,E)$. By induction hypothesis one can assume that $\sigma'$-stable $A$ with $\ext^1(A,A)<\ext^1(E,E)$
are $\sigma$-stable of the same phase. So the difficult case is the one that $E$ has a $\sigma'$-stable factor with the same $\ext^1$ and one needs to derive a contradiction here, somehow.
%
\end{remark}

\section{The spherical metric}\label{sect:Sphmetr}

\subsection{} We shall define a `spherical' version of Bridgeland's metric (see Section \ref{subsect:stabMet}). Instead of testing all object in
$\kd$ only spherical objects are  taken into account.

We start out with the space of slicings. As before, $\ks$ denotes the set of spherical objects in a
K3 category $\kd$.

\begin{definition}
For two slicings $\kp$ and $\kp'$ one defines 
$$f_\ks(\kp,\kp'):=\sup_{0\ne A\in\ks}\{|\phi^+(A)-{\phi'}^+(A)|,|\phi^-(A)-{\phi'}^-(A)|\}.$$
\end{definition}

Clearly, $f_\ks(\kp,\kp')\leq f(\kp,\kp')$ (see Section \ref{subsect:stabMet} for the definition of $f$) and thus the standard topology is a priori finer than the
one defined by $f_\ks$. 

\begin{remark}\label{rem:f=0}
i) Note that  $f_\ks$ for a general K3 category $\kd$  will usually not be a generalized metric and possibly  not
even well-defined. Eg.\ if $\kd$ has no or too few spherical objects, then $f_\ks$ is not defined (although one could set it constant zero
in this case) or one could have $f_\ks(\kp,\kp')=0$ without $\kp=\kp'$.

ii) Note that $f_\ks(\kp,\kp')=0$ if and only if $\kp(\phi)\cap\ks=\kp'(\phi)\cap\ks$ for all $\phi$. The `only if' is obvious. 
For the other direction, consider $A\in\ks$ with $\kp$-stable factors $A_1,\ldots,A_n$ of phase $\phi_1\geq\ldots\geq\phi_n$,
which are spherical  by Corollary \ref{cor:HMS1}. If
$\kp(\phi)\cap\ks=\kp'(\phi)\cap\ks$ for all $\phi$, then $A_i\in\kp'(\phi_i)$ and hence $\phi^\pm(A)={\phi'}^\pm(A)$.
This proves $f_\ks(\kp,\kp')=0$.

Note that  by Proposition \ref{prop:equivcond}
$f_\ks(\kp,\kp')=0$ is also equivalent to the condition
$\kp(\phi)^s\cap\ks=\kp'(\phi)^s\cap\ks$ for all $\phi$.  
\end{remark}

Consider two stability conditions $\sigma=(\kp,Z),\sigma'=(\kp',Z')$ on $\kd$.

\begin{definition}\label{def:sphmetr}
The \emph{spherical metric} $d_\ks(\sigma,\sigma')$ is defined 
as $$d_\ks(\sigma,\sigma'):=\max\{f_\ks(\kp,\kp'),|Z-Z'|\}.$$
\end{definition}

\subsection{} Consider a complex
projective K3 surface $X$ and let $\Sigma$ be Bridgeland's distinguished connected component
of the space $\Stab(X)$ of locally finite numerical stability conditions on $\Db(X)$ (see Section \ref{subsect:Stand}).
Let  $\sigma=(\kp,Z)$ and $\sigma'=(\kp',Z')$ be stability conditions contained in $\Sigma$.
The following is the analogue of \cite[Lem.\ 6.4]{BSt}.

\begin{prop}\label{prop:fdisc}
i) If $d_\ks(\sigma,\sigma')=0$, then $\sigma=\sigma'$.

ii) If  $Z=Z'$ and $f_\ks(\kp,\kp')<1$, then $\sigma=\sigma'$.
\end{prop}

\begin{proof}
Let us first prove that i) implies ii). Here we adapt the original arguments in \cite[Lem.\ 6.4]{BSt},
avoiding non-spherical objects.
So we have to show that  $Z=Z'$ and  $f_\ks(\kp,\kp')<1$
imply $f_\ks(\kp,\kp')=0$, i.e.\ that for all $\phi$ one has $\kp(\phi)\cap\ks=\kp'(\phi)\cap \ks$.
Assuming $A\in\kp(\phi)\cap\ks$, we have to prove $A\in\kp'(\phi)$. If $A\in\kp'(>\phi)$, then from
$f_\ks(\kp,\kp')<1$ one deduces $A\in\kp'(\phi,\phi+1)$. The latter would contradict $Z(A)=Z'(A)$.
Similarly one excludes the case $A\in\kp'(<\phi)$.

Consider the $\sigma'$-stable factors
 $A_1,\ldots,A_k$ of $A$ with $\phi'_1:=\phi'(A_1)\geq\ldots\geq\phi'_k:=\phi'(A_k)$.
Note that $A_1,\ldots,A_k$ are again spherical and
 that $\phi+1>\phi'_1\geq\ldots\geq\phi'_k>\phi-1$.
We have dealt already with the cases that $\phi>\phi'_1$ or $\phi_k'>\phi$. So we may assume
$\phi'_1\geq\phi\geq\phi'_k$ and have to show equality. If exactly one of the inequalities is strict, then $Z=Z'$ yields a contradiction.
So we may assume $\phi+1>\phi'_1\geq\ldots\geq\phi'_\ell>\phi\geq\phi'_{\ell+1}\geq\ldots\geq\phi'_k>\phi-1$
for some $1\leq\ell\leq k$.
 
 Breaking the filtration at the $\ell$-th step yields an  exact triangle $B_1\to A\to B_2$, i.e.\ the $\sigma'$-stable factors of $B_1$
 are $A_1,\ldots,A_\ell$ and the $\sigma'$-stable factors of
 $B_2$ are $A_{\ell+1},\ldots,A_k$.
 
 One now proves that $B_1\in\kp(>\phi-1)$ and $B_2\in\kp(\leq\phi+1)$. If $B_1\not\in\kp(>\phi-1)$,
 then there exists a $\sigma$-stable
 object $C$ of phase $\phi(C)\leq\phi-1$ such that $\Hom(B_1,C)\ne0$. Then for at least one of the $\sigma'$-stable
 factors $A_1,\ldots,A_\ell$ of $B_1$, say $A_{i_0}$, one has $\Hom(A_{i_0},C)\ne0$ and, in a next step,
 for one $\sigma$-stable factor $A_{i_0}'$ of $A_{i_0}$ one has $\Hom(A_{i_0}',C)\ne0$. Since the spherical $A_{i_0}$ is
 $\sigma'$-stable of phase $\phi'_{i_0}>\phi$, its $\sigma$-stable factor $A_{i_0}'$ has phase
 $\phi(A_{i_0}')\in(\phi'_{i}-1,\phi'_{i_0}+1)$. Hence $\phi(A_{i_0}')>\phi-1$, which contradicts $\Hom(A_{i_0}',C)\ne0$.
 The argument to prove $B_2\in\kp(\leq\phi+1)$ is similar.

Next one shows that $B_1\in\kp(>\phi-1)$ excludes $B_1\in\kp(\leq\phi)$. Indeed, otherwise $B_1\in\kp(\phi-1,\phi]$ and
hence $Z(B_1)\in\exp(i\pi\varphi)\RR_{>0}$ for some $\varphi\in(\phi-1,\phi]$. But $Z(B_1)=Z'(B_1)=\sum_{i_1}^\ell Z'(A_i)\in
\sum_{i=1}^\ell\exp(i\pi\phi_i')\RR_{>0}$ with $\phi_i'\in(\phi,\phi+1)$. Contradiction.

Now one concludes as in \cite{BSt}. As $B_1\not\in\kp(\leq\phi)$, there exists a $\sigma$-stable object
$C$ of phase $\phi(C)>\phi$ with $\Hom(C,B_1)=0$. Since $A\in\kp(\phi)$ one has $\Hom(C,A)=0$ and
hence $\Hom(C,B_2[-1])\ne0$.  But the latter is excluded due to $B_2[-1]\in\kp(\leq\phi)$.

%
%
%
To prove i) one observes that $d_\ks(\sigma,\sigma')=0$ implies $Z=Z'$ and $f_\ks(\kp,\kp')=0$ and that
by Remark \ref{rem:f=0} the latter is equivalent to $\kp(\phi)\cap \ks=\kp'(\phi)\cap \ks$ for all $\phi\in\RR$. From
Theorem \ref{thm:main} one concludes $\sigma=\sigma'$.
\end{proof}

\begin{cor}\label{cor:spmetr}
The classical metric $d$ and the spherical
metric $d_\ks$ define equivalent topologies on $\Sigma$.
\end{cor}

\begin{proof}
By Proposition \ref{prop:fdisc}, the projection $\Sigma\to \kp^+_0(X)$ is a local homeomorphism also for the
topology induced by $d_\ks$.
\end{proof}

\subsection{Stability conditions on spherical collections}\label{sect:bfks}
Ideally, we would like to talk about stability conditions on the set $\ks$, possibly viewed with its structure as
a $\CC$-linear category or with the binary operation $(A,B)\mapsto T_A(B)$ induced by spherical twists. However, there 
does not seem a way around the $\sigma$-stable filtrations and, although all the stable factors $A_i$ in (\ref{eqn:stabledomp}) for $E\in\ks$ are
 spherical (and for the spherical metric one only needs $A_1$ and $A_m$), the filtrations as such are not intrinsic to $\ks$.

So instead we consider $\ks^*\subset\Db(X)$, the smallest full triangulated subcategory containing $\ks$.
In other words $\ks$ generates $\ks^*$ without taking direct summands. For details see \cite{HuyBB}.
As noted there, $\ks^*$ is a triangulated category with a reasonably small Grothendieck group $K(\ks^*)=N(X)\subset
H^*(X,\ZZ)$ (assuming $\rho(X)\geq2$), but without bounded t-structures (we are working over $\CC$!).
As explained in \cite{HuyBB}, the category $\ks^*$ is expected to be $\Db(X)$ for K3 surfaces over $\bar\QQ$.

\begin{remark}
Note in passing that the numerical Grothendieck group of $\ks^*$, i.e.\ $K(\ks^*):=K(\ks^*)/_{\sim_\ks}$, equals $K(X)/_\sim=N(X)$,
for the Mukai pairing is non-degenerated on $N(X)$.
\end{remark}

So even passing to $\ks^*$ will not allow us to speak about stability conditions on $\ks$ or, rather,
on $\ks^*$. For this reason
we allow ourselves to adapt the original notion as follows. Let $\kt$ be a K3 category and $\ks\subset\kt$
a generating collection of spherical objects invariant under shift. For our purposes take $\kt=\ks^*$.

\begin{definition}\label{def:stabmodi}
A \emph{stability condition} $\sigma$ on $\ks$ with respect to $\ks\subset\kt$ consists of an
additive stability function $Z:K(\kt)/_{\sim}\to\CC$ and subsets $\ks(\phi)\subset\ks$, $\phi\in\RR$ satisfying the following conditions:
i)  $\ks(\phi)[1]=\ks(\phi+1)$, ii) $\Hom(\ks(\phi_1),\ks(\phi_2))=0$ for $\phi_1>\phi_2$, iii) $Z(\ks(\phi))\subset\exp(i\pi\phi)\RR_{>0}$, 
and  iv) for every $E\in\ks$ there exists a filtration as in (\ref{eqn:stabledomp}) with $A_i\in\ks(\phi_i)$.
\end{definition}

\begin{remark}
In order to think of this notion as a stability condition on $\ks$, i.e.\ independent of $\kt$, one would need
some kind of `formality' statement saying that $\kt$ is uniquely determined by the $\CC$-linear
category $\ks$. For certain `spherical configurations' this is indeed true (cf.\ \cite{KY,STh,Th}). In our
context one would in particular have to decide whether any $\CC$-linear equivalence $\ks_{X_1}\cong\ks_{X_2}$
for  the spherical collections $\ks_{X_i}\subset\Db(X_i)$, $i=1,2$, of two K3 surfaces $X_1,X_2$ always extends to
an exact equivalence $\Db(X_1)\cong\Db(X_2)$ (see \cite{HuyBB}).

Note that for a generic non-projective K3 surface $\ks$ consists of shifts of $\ko_X$ (see \cite{HMS}). In this
case, $\ks^*$ is the unique K3 category generated by a spherical object (cf.\ \cite{KY}).
\end{remark}

Let $\Stab(\ks):=\Stab(\ks\subset\kt)$ be the space of stability conditions on $\ks$ with respect to $\ks\subset\kt$ in the sense of Definition \ref{def:stabmodi}.
It can be equipped with a generalized metric $d_\ks$ as in Definition \ref{def:sphmetr}.
 We do not intend to develop the theory here fully, but most of the arguments in \cite{BSt} can be adapted. A good
 example is maybe Proposition \ref{prop:fdisc}, which works in this setting.

In any case, it is obvious that for $\kt=\ks^*\subset\Db(X)$ the restriction of a stability condition on $\Db(X)$
to $\ks^*$ yields a stability condition on $\ks$ (with respect to $\ks\subset\ks^*$) in the above sense.
The induced map $$\Stab(X)\to\Stab(\ks)$$ is continuous with respect to the corresponding metrics. The pull-back of the metric on
$\Stab(\ks)$ yields the spherical metric $d_\ks$ on $\Stab(X)$. The main result can thus be reformulated as
 
\begin{cor}\label{cor:stabsph}
On the distinguished component the restriction yields an embedding
$$\Sigma\,\,\hookrightarrow \Stab(\ks)$$
which identifies the natural metric on $\Stab(\ks)$ with the spherical metric
$d_\ks$ on $\Sigma$.\qqed
\end{cor}

\appendix{}
\section{Group of autoequivalences}

The following remarks are largely independent of the rest of the paper, but can be seen as
a motivation for the study of $\Stab(X)\to\Stab(\ks)$.

\subsection{} We shall first fix (or recall) some notations. As before, $X$ denotes a complex projective K3 surface
and $\Db(X)$  its bounded derived category of coherent sheaves. Due to the Global Torelli theorem
the group $\Aut(X)$ of automorphisms of $X$ can be identified with a subgroup of all
Hodge isometries of $H^2(X,\ZZ)$. We will write this
as $$\Aut(X)\,\hookrightarrow {\rm O}(H^2(X,\ZZ)).$$ The group of \emph{transcendental automorphisms}
$\TAut(X)$ is by definition the subgroup of $\Aut(X)$ consisting of all $f\in\Aut(X)$   for which the induced
action $f^*\in{\rm O}(H^2(X,\ZZ))$ is trivial on the algebraic part $\NS(X)$. Thus,
$$\TAut(X)\,\hookrightarrow{\rm O}(T(X)),$$ where $T(X)\subset H^2(X,\ZZ)$ is the transcendental lattice.
It is known that $\TAut(X)$ is a finite group.

The group of linear exact autoequivalences of $\Db(X)$ is denoted $\Aut(\Db(X))$. It comes with a natural representation
$$\rho:\Aut(\Db(X))\to {\rm O}(\widetilde H(X,\ZZ)),$$ which is not injective and we shall denote its kernel by
$\Aut_0(\Db(X))$. The description of the image of $\rho$ was completed in \cite{HMSDuke}. 
Since any autoequivalence naturally acts on $\Stab(X)$, one can define
the two  subgroups $$\Aut^\Sigma(\Db(X))\subset\Aut(\Db(X))\text{ and }
\Aut^\Sigma_0(\Db(X))\subset\Aut_0(\Db(X))$$
of autoequivalences that respect the distinguished component $\Sigma\subset\Stab(X)$. Conjecturally,
one has $\Sigma=\Stab(X)$ or, less optimistic, $\Aut^\Sigma(\Db(X))=\Aut(\Db(X))$.

\subsection{}
Instead of letting an automorphism of $X$ or an autoequivalence of $\Db(X)$ act on the cohomology, we can study its action on the collection of spherical 
objects $\ks\subset{\rm Ob}(\Db(X))$. We shall denote these \emph{spherical actions} by
$$\tau:\Aut(X)\to\Aut(\ks)\text{ and } \tau:\Aut(\Db(X))\to\Aut(\ks).$$

Note that the set $v(\ks)\subset \widetilde H(X,\ZZ)$ of  Mukai vectors
of all spherical objects generates the algebraic part $N(X)$ of $\widetilde H(X,\ZZ)$. This immediately
shows 
\begin{equation}\label{eqn:TAut}\ker(\tau:\Aut(X)\to\Aut(\ks))\subset\TAut(X).
\end{equation}

\begin{remark}
Presumably equality holds in (\ref{eqn:TAut}), but the only thing that seems obvious is the following. Let $f\in\TAut(X)$ and
suppose $A\in\ks$ is a spherical object that is stable with respect to some stability condition $\sigma\in\Sigma$.
Then $f^*A\cong A$.

For spherical sheaves which are $\mu$-stable with respect to some ample line bundle $H$ on $X$ 
this is due to a well-known argument of Mukai. If $A$ is $\mu_H$-stable, then $f^*A$ is $\mu_{f^*H}$-stable. Since $f^*H=H$ for a transcendental $f$,
both sheaves $A$ and $A'=f^*A$ are  $\mu_H$-stable. Moreover, they have the same Mukai vector
and hence $\chi(A,A')=2$ which shows that there must exist a non-trivial homomorphism
between $A$ and $A'$. The latter together with the stability of the two sheaves yields $A\cong A'=f^*A$.

For the general case of a spherical object $A$ that is stable with respect to some $\sigma\in\Sigma$, one
uses that $f^*\sigma=\sigma$ (see proof of Lemma \ref{lem:App2}) and argues as above.
One could try to deal with an arbitrary spherical object by applying the above to its stable factors (with respect to
some $\sigma\in\Sigma$), but due to the non-trivial action of $f$ on the Hom-spaces this is not obvious.
\end{remark}

\begin{lem}\label{lem:App1}
For the spherical representations one has
$$\ker(\tau:\Aut(\Db(X))\to\Aut(\ks))=\ker(\tau:\Aut(X)\to\Aut(\ks))\subset\TAut(X).$$
\end{lem}

\begin{proof} 
Let $\Phi\in\Aut(\Db(X))$ act trivially on $\ks$. In particular, $\Phi$ leaves invariant powers $L^i$ of an ample line bundle $L$. The action on the graded ring $\bigoplus H^0(X,L^i)$ is induced by an automorphism
$f\in\Aut(X,L)$. Hence $\Phi$ and $f^*$ define two autoequivalences which are isomorphic
on the full subcategory given by the ample sequence $\{L^i\}$. By a result of Bondal and Orlov (see \cite{BO} or \cite[Prop.\ 4.23]{HFM}), this immediately yields $\Phi=f^*$. But then $f\in\TAut(X)$.
\end{proof}


\subsection{}
The groups $\Aut(X)$ and $\Aut(\Db(X))$ both act on $\Stab(X)$. 
We denote this action by $$\kappa:\Aut(\Db(X))\to\Aut(\Stab(X)).$$ The main result of
\cite{BK3} says that  the subgroup
$\Aut_0^\Sigma(\Db(X))$ is via $\kappa$ identified with the group of deck transformations
of $\Sigma\to\kp^+_0(X)$:
$$\kappa:\Aut_0^\Sigma(X)\congpf{\rm Gal}(\Sigma/\kp_0^+(X)).$$

\begin{lem}\label{lem:App2}
For the action on the space of stability conditions one
has $$\ker(\kappa:\Aut^\Sigma(\Db(X))\to\Aut(\Sigma))=\ker(\kappa:\Aut(X)\to\Aut(\Sigma))=\TAut(X).$$
\end{lem}

\begin{proof} Let first $f\in\Aut(X)$ and consider a standard stability condition $\sigma$ 
with stability function $Z(E)=\langle \exp(B+i\omega),v(E)\rangle$.
Then $f^*\sigma$ is a stability function for which by definition all $f^*k(x)$ are stable.
Hence all point sheaves $k(y)$ are again stable and for $f\in\TAut(X)$ the stability function remains
unchanged under pull-back. Thus, for standard stability conditions $\sigma$ one has $f^*\sigma=\sigma$.
In particular, $f^*$ preserves the distinguished component $\Sigma$ and acts on it as a deck-transformation
with fixed points. Hence $f^*=\id$ on $\Sigma$.

Conversely, if $\Phi$ acts trivially on $\Sigma$, then $\Phi$ acts trivially on $\NS(X)$.
By the Global Torelli theorem the induced action on $T(X)$ is of the form $f^*$ for some $f\in\TAut(X)$. Changing
$\Phi$ by the inverse of $f^*$, we may assume that $\Phi$ acts trivially on $\widetilde H(X,\ZZ)$. But then
$\Phi\in\ker(\kappa:\Aut_0^\Sigma(\Db(X))\to{\rm Gal}(\Sigma/\kp^+_0(X))$ and thus $\Phi=\id$. 
\end{proof}

The observation that the kernels
of the two actions  
$$\tau:\Aut(\Db(X))\to\Aut(\ks)) ~\text{  and } ~\kappa:\Aut(\Db(X))\to
\Aut(\Sigma)$$ essentially
coincide hints at the deeper  that stability conditions in $\Sigma$ are determined by their behavior
with respect to $\ks$ which is expressed by Theorem \ref{thm:main} and Corollary \ref{cor:stabsph}.




{\footnotesize \begin{thebibliography}{mm}


\bibitem{BO} A.\ Bondal,  D.\ Orlov, 
\em Reconstruction of a variety from the derived category and groups of autoequivalences. \em
Comp.\  Math.\ 125 (2001), 327--344. 


\bibitem{Brav} C.\ Brav, H.\ Thomas
\em Braid groups and Kleinian singularities. \em
 arXiv:0910.2521.

\bibitem{BSt} T.\ Bridgeland
\em Stability conditions on triangulated categories. \em
Ann.\ of Math.\ 166 (2007),  317--345.

\bibitem{BK3} T.\ Bridgeland
\em Stability conditions on K3 surfaces. \em Duke Math.\ J.\ 141 (2008),  241--291. 

\bibitem{BKl} T.\ Bridgeland
\em Stability conditions and Kleinian singularities. \em
Int.\ Math.\ Res.\ Not.\ 21 (2009), 4142--4157. 


\bibitem{HFM} D.\ Huybrechts
\em Fourier--Mukai transforms in algebraic geometry. \em Oxford
Mathematical Monographs (2006).

\bibitem{HMS} D.\ Huybrechts, E.\ Macr\`i, P.\ Stellari
\em Stability conditions for generic $K3$ categories. \em Compos.\ Math.\ 144 (2008), 134--162.

\bibitem{HMSDuke} D.\ Huybrechts, E.\ Macr\`i, P.\ Stellari
\em Derived equivalences of K3 surfaces and orientation. \em
Duke Math.\ J.\ 149 (2009), 461--507. 

\bibitem{HCH} D.\ Huybrechts
\em Chow groups and derived categories of K3 surfaces. \em
J.\ EMS.\ 12 (2010), 1533--1551.

\bibitem{HuyBB} D.\ Huybrechts
\em A note on the Bloch--Beilinson conjecture for K3 surfaces and spherical objects. \em
Preprint (2010).


\bibitem{Ishi} A.\ Ishii, K.\ Ueda, H.\ Uehara
\em Stability conditions on $A_n$-singularities. \em
J.\ Diff.\ Geom.\ 84 (2010), 87--126. 

\bibitem{KY}
B.\ Keller, D.\ Yang, Dong, G.\ Zhou
\em The Hall algebra of a spherical object. \em
J.\ Lond.\ Math.\ Soc.\ (2) 80 (2009),  771--784.

\bibitem{Mu} S.\ Mukai
\em On the moduli space of bundles on K3 surfaces, I. \em
 In: Vector Bundles on Algebraic Varieties, Bombay
(1984), 341-413.


\bibitem{STh} P.\ Seidel, R.\ Thomas
\em Braid group actions on derived categories of coherent sheaves. \em
Duke Math.\ J.\ 108 (2001), 37--108.

\bibitem{Th} R.\ Thomas
\em Stability conditions and the braid group. \em
Comm.\ Anal.\ Geom.\ 14 (2006),  135--161. 

\end{thebibliography}}
\end{document}